\documentclass[14pt]{article}
\usepackage{mathrsfs}
\usepackage{amsthm}
\usepackage{amsmath}
\usepackage{graphicx}
\usepackage{color}
\usepackage{amsfonts}
\usepackage[text={140mm,210mm},left=45mm,vmarginratio=1:1]{geometry}
\newtheorem{theorem}{Theorem}[section]

\newtheorem{lemma}[theorem]{Lemma}

\newtheorem{corollary}[theorem]{Corollary}
\newtheorem{question}[theorem]{Question}
\numberwithin{equation}{section}
\normalsize

\begin{document}
\title{\textbf{Critical value for the contact process with random edge weights on regular tree}}

\author{Xiaofeng Xue \thanks{\textbf{E-mail}: xuexiaofeng@ucas.ac.cn \textbf{Address}: School of Mathematical Sciences, University of Chinese Academy of Sciences, Beijing 100049, China.}\\ University of Chinese Academy of Sciences}

\date{}
\maketitle

\noindent {\bf Abstract:}

In this paper we are concerned with contact processes with random
edge weights on rooted regular trees. We assign i.i.d weights on
each edge on the tree and assume that an infected vertex infects its
healthy neighbor at rate proportional to the weight on the edge
connecting them. Under the annealed measure, we define the critical
value $\lambda_c$ as the maximum of the infection rate with which
the process will die out and define $\lambda_e$ as the maximum of
the infection rate with which the process dies out at exponential
rate. We show that these two critical values satisfy an identical
limit theorem and give an precise lower bound of $\lambda_e$. We
also study the critical value under the quenched measure. We show
that this critical value equals that under the annealed measure or
infinity according to a dichotomy criterion. The contact process on
a Galton-Watson tree with binomial offspring distribution is a
special case of our model.


\section{Introduction}\label{section one}
In this paper, we are concerned with contact processes with random
edge weights on regular trees. For each integer $N\geq 1$, we denote
by $\mathbb{T}^N$ the rooted regular tree where the root $O$ has
degree $N$ and other vertices have degree $N+1$. That is to say,
each vertex produces $N$ children and the root $O$ has no ancestor
while each other vertex has a father. The following picture
describes a local area of $\mathbb{T}^4$.

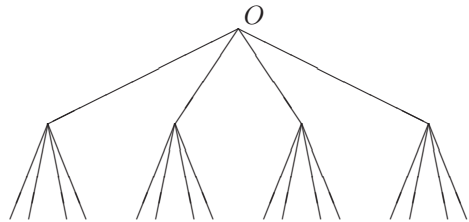
\begin{center}
\begin{picture}(250,110)
\put(112,77){$O$} \put(110,75){\line(2,-1){72}}
\put(110,75){\line(-2,-1){72}} \put(110,75){\line(2,-3){24}}
\put(110,75){\line(-2,-3){24}} \put(182,39){\line(1,-5){7.2}}
\put(182,39){\line(2,-5){14.4}} \put(182,39){\line(-2,-5){14.4}}
\put(182,39){\line(-1,-5){7.2}} \put(134,39){\line(1,-5){7.2}}
\put(134,39){\line(2,-5){14.4}} \put(134,39){\line(-2,-5){14.4}}
\put(134,39){\line(-1,-5){7.2}} \put(86,39){\line(1,-5){7.2}}
\put(86,39){\line(2,-5){14.4}} \put(86,39){\line(-2,-5){14.4}}
\put(86,39){\line(-1,-5){7.2}} \put(38,39){\line(1,-5){7.2}}
\put(38,39){\line(2,-5){14.4}} \put(38,39){\line(-2,-5){14.4}}
\put(38,39){\line(-1,-5){7.2}}
\end{picture}
\end{center}

For any vertices $x,y\in \mathbb{T}^N$, we denote by $x\sim y$ when
there is an edge connecting them. We denote by $\mathbb{E}^N$ the
set of edges on $\mathbb{T}^N$.

Let $\rho$ be a non-negative random variable such that $P(\rho\leq
M)=1$ for some $M\in (0,+\infty)$ and $P(\rho>0)>0$.
$\{\rho(e)\}_{e\in \mathbb{E}^N}$ are i. i. d. random variables such
that for each $e\in \mathbb{E}^N$, $\rho(e)$ and $\rho$ have the
same probability distribution. For $e\in \mathbb{E}^N$ with
endpoints $x,y\in \mathbb{T}^N$, we write $\rho(e)$ as $\rho(x,y)$.
When $\{\rho(e)\}_{e\in \mathbb{E}^N}$ is given, the contact process
with edge weights $\{\rho(e)\}_{e\in \mathbb{E}^N}$ is a spin system
with state space $\{0,1\}^{\mathbb{T}^N}$ and flip rates function
given by
\begin{equation}\label{equ 1.1 flip rate}
c(x,\eta)=
\begin{cases}
1 &\text{~if~}\eta(x)=1,\\
\lambda\sum_{y:y\sim x}\rho(x,y)\eta(y)&\text{~if~}\eta(x)=0
\end{cases}
\end{equation}
for any $(x,\eta)\in \mathbb{T}^N\times \{0,1\}^{\mathbb{T}^N}$,
where $\lambda$ is a positive parameter called the infection rate.

The assumption $P(\rho\leq M)=1$ ensures the existence of our
process according to the basis theory constructed in \cite{Har1972}
and \cite{Lig1972}.

Intuitively, the process describes the spread of an infection
disease. Vertices in state $1$ are infected individuals while
vertices in state $0$ are healthy. An infected individual waits for
an exponential time with rate $1$ to recover. A healthy vertex $x$
is infected by its infected neighbor $y$ at a rate proportional to
$\rho(x,y)$. That is to say, the larger $\rho(x,y)$ is, the faster
the disease spreads from $y$ to $x$.

When $\rho\equiv 1$, our model degenerates to the classic contact
process, which is introduced in \cite{Har1974} by Harris. In
\cite{Pem1992}, Pemantle first considers contact processes on trees.
The two books \cite{Lig1985} and \cite{Lig1999} written by Liggett
give a detailed introduction for the study of classic contact
processes on lattices and trees.

When $P(\rho=1)=1-P(\rho=0)=p\in (0,1)$, then our model turns into
the contact process on a Galton-watson tree with binomial offspring
distribution $B(N,p)$ and also can be seen as contact process on
open clusters of bond percolation on tree. In \cite{Pem2001},
Pemantle and Stacey study contact processes and branching random
walks on Galton-Watson trees. They show that on some Galton-Watson
trees the branching random walk has one phase transition while the
contact process has two. Contact processes on clusters of bond
percolation on lattices are studied by Chen and Yao in
\cite{Chen2009}. They show that the complete convergence theorem
holds.

In this paper, we are concerned with contact processes with random
edge weights. It is also interesting to consider the process with
random vertex weights. In detail, each vertex $x$ is assigned a
weight $\rho(x)$. Infected vertex $x$ infects healthy neighbor $y$
at rate proportional to $\rho(x)\rho(y)$. This model concludes
contact process on clusters of site percolation as a special case.
In \cite{Ber2011}, Bertacchi, Lanchier and Zucca study contact
processes on $C_{\infty}\times K_N$, where $C_\infty$ is the
infinite open cluster of site percolation and $K_N$ is a complete
graph with $N$ vertices. Criterions to judge whether the process
will survive are given in \cite{Ber2011}. Contact processes with
random vertex weights on complete graphs are introduced in
\cite{Pet2011} by Peterson. In \cite{Pet2011}, it is shown that the
critical value of the model is inversely proportional to the second
moment of the vertex weight. Xue extends this result to the case
where the graph is oriented lattice in \cite{Xue2015}. In
\cite{Xue2013}, Xue studies contact processes with random vertex
weights on general regular graphs and obtains a lower bound of the
critical value of the model.

\section{Main results}\label{section two}
In this section we give main results of this paper. First we
introduce some notations and definitions. For each $N\geq 1$, we
assume that $\{\rho(e)\}_{e\in \mathbb{E}^N}$ are defined on the
probability space $\big(\Omega_N, \mathcal{F}_N, \mu_N\big)$. We
write $(\Omega_N,\mathcal{F}_N,\mu_N)$ briefly as
$(\Omega,\mathcal{F},\mu)$ when there is no misunderstanding. For
any $\omega\in \Omega$, we denote by $P_\lambda^{\omega}$ the
probability measure of the contact process on $T^N$ with edge
weights $\{\rho(e,\omega)\}_{e\in \mathbb{E}^N}$ and infection rate
$\lambda$. $P^{\omega}_{\lambda}$ is called the quenched measure.
The expectation operator with respect to $P^{\omega}_\lambda$ is
denoted by $E^{\omega}_\lambda$. We define
\[
P_\lambda^N(\cdot)=\int P^{\omega}_\lambda(\cdot)\mu_N(d\omega),
\]
which is called the annealed measure. The expectation operator with
respect to $P_\lambda^N$ is denoted by $E^N_\lambda$. When there is
no misunderstanding, we write $P_\lambda^N$ and $E_\lambda^N$
briefly as $P_\lambda$ and $E_\lambda$.

For any $t\geq 0$, we denote by $\eta_t$ the configuration of the
contact process at moment $t$. The value of vertex $x$ at moment $t$
is denoted by $\eta_t(x)$. For any $t>0$, let
\[
C_t=\{x\in T^N:\eta_t(x)=1\}
\]
be the set of infected vertices at $t$. We write $C_t$ as $C_t^O$
when $C_0=\{O\}$.

Since $\emptyset$ is an absorbed state of the process
$\{C_t\}_{t\geq 0}$ and the contact process is an attractive spin
system (see section 3.2 of \cite{Lig1985}), for any
$\lambda_1>\lambda_2$,
\begin{equation}\label{equ 2.1}
P_{\lambda_1}(\forall~t\geq 0,C_t^O\neq \emptyset)\geq
P_{\lambda_2}(\forall~t\geq 0,C_t^O\neq \emptyset).
\end{equation}
By \eqref{equ 2.1}, it is reasonable to define the following
critical value. For each $N\geq 1$, we define
\begin{equation}\label{equ 2.2 anneal critical value}
\lambda_c(N)=\sup\{\lambda:P_\lambda^N(\forall~t\geq 0,C_t^O\neq
\emptyset)=0\}.
\end{equation}
We write $\lambda_c(N)$ as $\lambda_c$ when there is no
misunderstanding.

Supposing that only $O$ is infected at $t=0$, then when
$\lambda<\lambda_c$, with probability one there will be no infected
vertices eventually, which means that the disease dies out. When
$\lambda>\lambda_c$, with positive probability there will be always
some vertices in the infected state, which means that the disease
survives. The case of $\lambda=\lambda_c$ is difficult. In
\cite{Bez1990}, Bezuidenhout and Grimmett show that the critical
classic contact process on lattice dies out. We guess same
conclusion holds for our model but have not find a way to prove it
yet.

When $\lambda<\lambda_c$,
\[
\lim_{t\rightarrow+\infty}P_\lambda(C_t^O\neq \emptyset)=0.
\]
It is natural to ask whether $P_\lambda(C_t^O\neq \emptyset)$
converges to $0$ at an exponential rate. So it is natural to define
the following critical value. For any $N\geq1$, we define
\begin{equation}\label{equ 2.3 anneal exponential critical value}
\lambda_e(N)=\sup\{\lambda:\limsup_{t\rightarrow+\infty}\frac{1}{t}\log
P_\lambda^N(C_t^O\neq \emptyset)<0\}.
\end{equation}
It is obviously that $\lambda_e\leq \lambda_c$. Does $\lambda_e$
equal $\lambda_c$? Section 6.3 of \cite{Lig1985} shows that the
answer is positive for classic contact process on $Z$. We have no
idea whether $\lambda_e=\lambda_c$ for our model.

Now we give our main results. Our first result is a criterion to
judge whether $\lambda_c\in(0,+\infty)$.

\begin{theorem}\label{theorem 2.0}
If $P(\rho>0)=1$, then for each $N\geq2$, $0<\lambda_c(N)<+\infty
 $. If $P(\rho>0)<1$, then $0<\lambda_c(N)<+\infty$
 for $N>1/P(\rho>0)$ and $\lambda_c(N)=+\infty$ for $N\leq
 1/P(\rho>0)$.
\end{theorem}

We can not judge whether $\lambda_c<+\infty$ for the case where
$N=1$ and $P(\rho>0)=1$. We guess in this case there is no common
conclusion. More information about the distribution of $\rho$ is
needed. For example, if there exists $\epsilon>0$ such that
$P(\rho>\epsilon)=1$, then it is easy to see that $\lambda_c\in
(0,+\infty)$ since classic contact process on $Z$ has finite
critical value (see Section 6.1 of \cite{Lig1985} and
\cite{Lig1995}).

To describe $\lambda_c$ and $\lambda_e$ more accurately, we obtain a
limit theorem of $\lambda_c$, $\lambda_e$ and a precise lower bound
of $\lambda_e$.

\begin{theorem}\label{theorem 2.1 limit theorem for two critical}
For $\rho$ satisfies that $P(\rho>0)>0$ and $P(0\leq \rho\leq M)=1$
for some $M\in (0,+\infty)$,
\begin{equation}\label{equ 2.4 limit theorem for two critical}
\lim_{N\rightarrow+\infty}N\lambda_c(N)=\lim_{N\rightarrow+\infty}N\lambda_e(N)=\frac{1}{E\rho}.
\end{equation}
Furthermore,
\begin{equation}\label{equ 2.5 lower bound of lambda-e}
\lambda_e(N)\geq\big(NE\rho+\frac{M^2}{E\rho}\big)^{-1}.
\end{equation}
\end{theorem}

Theorem \ref{theorem 2.1 limit theorem for two critical} show that
$\lambda_c,\lambda_e\approx 1/(NE\rho)$, which is inversely
proportional to the degree of the root and the mean of the edge
weight.

Let us see some examples. When $\rho\equiv 1$, Theorem \ref{theorem
2.1 limit theorem for two critical} shows that
\[
\lim_{N\rightarrow+\infty}N\lambda_c(N)=1
\]
and $\lambda_c(N)\geq 1/(N+1)$, which is the estimation of critical
value for classic contact process on regular tree given in
\cite{Pem1992}.

When $P(\rho=1)=1-P(\rho=0)=p\in (0,1)$, Theorem \ref{theorem 2.1
limit theorem for two critical} gives the estimation of critical
value for contact processes on Galton-Watson tree with binomial
offspring distribution $B(N,p)$ that
\[
\lim_{N\rightarrow+\infty}Np\lambda_c(N)=1
\]
and
\[
\lambda_c(N)\geq\frac{1}{Np+1/p}.
\]
These two estimations do not occur in former references.

The critical value $\lambda_c$ is defined under the annealed
measure. It is natural to consider the critical value of the process
with fixed edge weights $\{\rho(e,\omega)\}_{e\in \mathbb{E}^N}$ for
some $\omega\in \Omega$. Hence, for any $\omega\in \Omega_N$, we
define
\begin{equation}\label{equ 2.6 quench critical value}
\widehat{\lambda}_c(\omega,N)=\sup\{\lambda:P_\lambda^\omega(\forall~t,
C_t^{O}\neq\emptyset)=0\}.
\end{equation}

For $\omega\in \Omega_N$, if there is a cut-off $\Pi$ of
$\mathbb{T}^N$ separating $O$ from infinity such that
$\rho(e,\omega)=0$ for each $e\in \Pi$, then it is easy to see that
$\widehat{\lambda}_c(\omega,N)=+\infty$. We can show that except
this case, $\widehat{\lambda}_c(\omega,N)=\lambda_c(N)$, which means
the critical values under the annealed measure and quenched measure
are equal. To introduce our result rigorously, we introduce some
notations and definitions.

For any $\omega\in \Omega_N$, we define
\begin{equation}\label{equ 2.7}
L(\omega)=\{e\in \mathbb{E}^N:\rho(e,\omega)>0\}.
\end{equation}
For each $x\in \mathbb{T}^N$, there is an unique path $p(O,x)$ from
$O$ to $x$ which does not backtrack. We write $O\rightarrow^\omega
x$ when and only when each edge of $p(O,x)$ belongs to $L(\omega)$.
We define
\begin{equation}\label{equ 2.7 two}
D(\omega)=\{x\in T^N:O\rightarrow^\omega x\}
\end{equation}
and
\begin{equation}\label{equ 2.8}
A_N=\{\omega:|D(\omega)|=+\infty\}.
\end{equation}
It is obviously that $D(\omega)$ forms a Galton-Watson tree with
offspring distribution $B(N,q)$ and $1-\mu_N(A_N)$ is the extinction
probability of the tree, where
\[
q=P(\rho>0)
\]
and
\[
P(B(N,q)=k)={N\choose k}q^k(1-q)^{N-k}
\]
for $1\leq k\leq N$.

Now we can give our result of the critical value under the quenched
measure.
\begin{theorem}\label{theorem 2.2 two critical values are equal}
If $P(\rho>0)=1$, then for each $N\geq 2$, there exists $K_N\in F_N$
such that $\mu_N(K_N)=1$ and
\[
\widehat{\lambda}_c(\omega,N)=\lambda_c(N)\in (0,+\infty)
\]
for any $\omega\in K_N$, where $\lambda_c(N)$ is the same as that in
\eqref{equ 2.2 anneal critical value}.

If $P(\rho>0)<1$, then when $N\leq 1/P(\rho>0)$,
\[
\widehat{\lambda}_c(\omega,N)=+\infty
\]
for any $\omega\in \Omega_N$. When $N>1/P(\rho>0)$, then
$\mu_N(A_N)>0$ and there exists $K_N\subseteq A_N$ such that
$\mu_N(A_N\setminus K_N)=0$ and
\[
\widehat{\lambda}_c(\omega,N)=\lambda_c(N)\in (0,+\infty)
\]
for any $\omega\in K_N$. For any $\omega\not\in A_N$,
\[
\widehat{\lambda}_c(\omega,N)=+\infty.
\]
\end{theorem}
In conclusion, theorem \ref{theorem 2.2 two critical values are
equal} shows that $\widehat{\lambda}_c(\omega,N)\in
\{\lambda_c(N),+\infty\}$ with probability one. Furthermore,
\[\{\omega:\widehat{\lambda}_c(\omega,N)=\lambda_c(N)\}=A_N\] and
\[\{\omega:\widehat{\lambda}_c(\omega,N)=+\infty\}=\Omega_N\setminus A_N\]
in the sense of ignoring a set with probability $0$.

The proofs of our main results are divided into three sections. In
Section \ref{section three}, we will give an upper bound of
$\lambda_c$, which shows that
$\limsup_{N\rightarrow+\infty}N\lambda_c(N)\leq 1/E\rho$. The core
idea is to compare the contact process with a SIR epidemic model.
This section also gives most part of the proof of Theorem
\ref{theorem 2.0} except showing that $\lambda_c>0$. In Section
\ref{section four}, we will prove that $\lambda_e(N)\geq
\big(NE\rho+\frac{M^2}{E\rho}\big)^{-1}$ and hence $\lambda_c>0$,
which accomplishes the proof of Theorem \ref{theorem 2.1 limit
theorem for two critical} and Theorem \ref{theorem 2.0}. The main
approach is to compare the contact process with the binary contact
path process introduced in \cite{Grif1983} by Griffeath. In
technique, we need to estimate the number of paths (may backtrack)
from $O$ with given length on the tree. We relate this problem to
simple random walk on regular tree. In Section \ref{section five},
we give the proof of Theorem \ref{theorem 2.2 two critical values
are equal}. Our approach is inspired by the classic method of
proving extinction criterion for Galton-Watson trees.

\section{Upper bound for $\lambda_c$}\label{section three}
In this section we will prove that
$\limsup_{N\rightarrow+\infty}N\lambda_c(N)\leq 1/E\rho$. The
following lemma gives an upper bound of $\lambda_c(N)$, which is
crucial for our proof.

\begin{lemma}\label{lemma 3.1}
If $\lambda$ satisfies that
\[
NE\big[\frac{\lambda \rho}{1+\lambda \rho}\big]>1,
\]
then
\[
\lambda_c(N)\leq \lambda.
\]
\end{lemma}
We give the proof of Lemma \ref{lemma 3.1} at the end of this
section. First we utilize Lemma \ref{lemma 3.1} to prove that
$\limsup_{N\rightarrow+\infty}N\lambda_c(N)\leq 1/E\rho$.

\proof[Proof of
$~\limsup\limits_{N\rightarrow+\infty}N\lambda_c(N)\leq 1/E\rho$]

For $\gamma>1$, let $\lambda=\frac{\gamma}{NE\rho}$, then
\[
NE\big[\frac{\lambda\rho}{1+\lambda\rho}\big]=\frac{\gamma}{E\rho}E\big[\frac{\rho}{1+\frac{\gamma\rho}{NE\rho}}\big].
\]
According to Domination Convergence Theorem,
\[
\lim_{N\rightarrow+\infty}\frac{\gamma}{E\rho}E\big[\frac{\rho}{1+\frac{\gamma\rho}{NE\rho}}\big]=\frac{\gamma}{E\rho}E\rho=\gamma>1.
\]
Therefore, for sufficiently large $N$ and
$\lambda=\frac{\gamma}{NE\rho}$,
\[
NE\big[\frac{\lambda\rho}{1+\lambda\rho}\big]>1.
\]
Therefore, according to Lemma \ref{lemma 3.1},
\[
\lambda_c(N)\leq \frac{\gamma}{NE\rho}
\]
for sufficiently large $N$ and hence
\[
\limsup_{N\rightarrow+\infty}N\lambda_c(N)\leq \frac{\gamma}{E\rho}.
\]
Since $\gamma$ is arbitrary, let $\gamma\rightarrow 1$ and the proof
is complete.

\qed

For the special case $P(\rho=1)=1-P(\rho=0)=p\in (0,1]$ and $N>1/p$,
Lemma \ref{lemma 3.1} gives a precise upped bound of $\lambda_c(N)$
that
\[
\lambda_c(N)\leq \frac{1}{Np-1},
\]
since $NE\big[\frac{\lambda\rho}{1+\lambda\rho}\big]=\frac{\lambda
Np}{1+\lambda}$.

According to Lemma \ref{lemma 3.1}, we can also judge whether
$\lambda_c<+\infty$.

\begin{corollary}\label{corollary 3.2}
If $P(\rho>0)=1$, then $\lambda_c(N)<+\infty$ for each $N\geq 2$. If
$P(\rho>0)<1$, then $\lambda_c(N)<+\infty$ for $N>1/P(\rho>0)$ and
$\lambda_c(N)=+\infty$ for $N\leq 1/P(\rho>0)$.
\end{corollary}

\proof

According to Domination Convergence Theorem,
\[
\lim_{\lambda\rightarrow+\infty}E\big[\frac{\lambda\rho}{1+\lambda\rho}\big]=P(\rho>0).
\]
Therefore, in the case where $P(\rho>0)=1$ and $N\geq 2$ and the
case where $P(\rho>0)<1$ and $N> 1/P(\rho>0)$,
\[
\lim_{\lambda\rightarrow+\infty}NE\big[\frac{\lambda\rho}{1+\lambda\rho}\big]>1
\]
and hence
\[
\lambda_c(N)<\lambda
\]
for sufficiently large $\lambda$ according to Lemma \ref{lemma 3.1}.
As a result, in these two cases,
\[
\lambda_c<+\infty.
\]
For the case where $P(\rho>0)<1$ and $N\leq 1/P(\rho>0)$, the
Galton-Watson tree with offspring distribution $B(N,P(\rho>0))$ is
extinct with probability one, since the mean of the number of sons
is at most one. As a result, $D(\omega)$ is finite with probability
one and the Markov process $\{C_t^{O}\}_{t\geq 0}$ is with finite
state space $\{A:A\subseteq D(\omega)\}$. Since $\emptyset$ is the
unique absorption state for $\{C_t^{O}\}_{t\geq 0}$, the process
will be frozen in state $\emptyset$ eventually. As a result, for any
$\lambda>0$,
\[
P_\lambda^{\omega}(\forall ~t\geq 0, C_t^O\neq \emptyset)=0
\]
for any $\omega\in\Omega$ except a set with probability zero and
hence
\[
P_\lambda^N(\forall ~t\geq 0, C_t^O\neq \emptyset)=0.
\]
Therefore, $\lambda_c>\lambda$ for any $\lambda>0$ and hence
\[
\lambda_c=+\infty.
\]

\qed

At last we give the proof of Lemma \ref{lemma 3.1}.

\proof[Proof of Lemma \ref{lemma 3.1}]

To control the size of $C_t$ from below, we introduce the following
SIR epidemic model with random edge weights. Let $\{\xi_t\}_{t\geq
0}$ be Markov process with state space $\{-1,0,1\}^{\mathbb{T}^N}$.
At $t=0$, $\xi_0(O)=1$ and $\xi_0(x)=0$ for each other $x\in
\mathbb{T}^N$. For any $t\geq0$, we define
\begin{align*}
&I_t=\{x\in \mathbb{T}^N:\xi_t(x)=1\}, S_t=\{x\in \mathbb{T}^N:\xi_t(x)=0\},\\
&R_t=\{x\in \mathbb{T}^N:\xi_t(x)=-1\}.
\end{align*}
Now we can identify $\xi_t$ with $(S_t,I_t,R_t)$. After the edge
weights $\{\rho(e)\}_{e\in \mathbb{E}^N}$ is given,
$\{(S_t,I_t,R_t)\}_{t\geq 0}$ evolves as follows. For each $x\in
I_t$, $(S_t,I_t,R_t)$ flips to $(S_t,I_t\setminus \{x\},R_t\cup
\{x\})$ with rate $1$. For any $x,y$ satisfy that $y$ is a son of
$x$, $x\in I_t$ and $y\in S_t$, $(S_t,I_t,R_t)$ flips to
$(S_t\setminus \{y\}, I_t\cup \{y\},R_t)$ at rate
$\lambda\rho(x,y)$.

Intuitively, $1,0,-1$ represent `infected', `healthy' and `removed'
respectively. An infected vertex waits for an exponential time with
rate one to become removed. A healthy vertex $y$ may be infected
when and only when its father $x$ is infected. $x$ infects $y$ at
rate proportional to $\rho(x,y)$. A removed vertex will stay in the
this state forever.

For $\{C_t\}_{t\geq 0}$, an infected vertex can infect any healthy
neighbor while for $\{\xi_t\}_{t\geq 0}$, an infected vertex can
only infect its sons. For $\{C_t\}_{t\geq 0}$, when an infected
vertex become healthy, it may be infected again while for
$\{\xi_t\}_{t\geq 0}$, when an infected vertex becomes removed, it
will never be infected again. As a result, according to the approach
of basic coupling (see section 2.1 of \cite{Lig1985}), it is easy to
see that
\[
I_t\subseteq C_t^{O}
\]
for any $t>0$ in the sense of coupling when the two processes with
same infection rate $\lambda$ and edge weights $\{\rho(e)\}_{e\in
\mathbb{E}^N}$. Therefore,
\begin{equation*}
P_\lambda^\omega(\forall~t,C_t^{O}\neq \emptyset)\geq
P_\lambda^\omega(\forall~t,I_t\neq \emptyset)
\end{equation*}
for any $\omega\in \Omega$ and hence
\begin{equation}\label{equ 3.1}
P_\lambda^N(\forall~t,C_t^{O}\neq \emptyset)\geq
P_\lambda^N(\forall~t,I_t\neq \emptyset)
\end{equation}
for any $t>0$ and $N\geq 1$.

We define
\[
I_{+\infty}=\bigcup_{t\geq 0}I_t
\]
as the set of vertices which have been infected. $I_t\neq \emptyset$
for any $t\geq 0$ if and only if there are infinite many vertices
which have been infected. Therefore,
\begin{equation}\label{equ 3.2}
\{\forall~t\geq 0,I_t\neq \emptyset\}=\{|I_{+\infty}|=+\infty\}.
\end{equation}
By \eqref{equ 3.1} and \eqref{equ 3.2},
\begin{equation}\label{equ 3.3}
P_\lambda^N(\forall~t\geq 0,C_t\neq \emptyset)\geq
P_\lambda^N(|I_{+\infty}|=+\infty).
\end{equation}
For $x\in T^N$ and a son $y$ of $x$, let $T_1$ be an exponential
time with rate $\lambda\rho(x,y)$ and $T_2$ be an exponential time
with rate $1$ and independent of $T_1$, then conditioned on $x$ is
infected, the probability that $x$ infects $y$ equals
\[
P(T_1<T_2)=\frac{\lambda\rho(x,y)}{1+\lambda\rho(x,y)}.
\]
As a result, under the annealed measure $P_\lambda^N$, the mean of
the number of infected sons of an infected vertex equals
\[
NE\big[\frac{\lambda\rho}{1+\lambda\rho}\big].
\]
As a result, under the annealed measure $P_\lambda^N$, $I_{+\infty}$
forms a Galton-Watson tree with an offspring distribution with mean
$ NE\big[\frac{\lambda\rho}{1+\lambda\rho}\big] $. According to the
extinction criterion of Galton-Watson trees,
\begin{equation}\label{equ 3.4}
P_\lambda^N(|I_{+\infty}|=+\infty)>0
\end{equation}
when $\lambda$ satisfies
$NE\big[\frac{\lambda\rho}{1+\lambda\rho}\big]>1$. Lemma \ref{lemma
3.1} follows from this fact and \eqref{equ 3.3}.

\qed

\section{Lower bound for $\lambda_e$}\label{section four}
In this section we will give lower bound of $\lambda_e$. First, we
give a lemma about simple random walk on $\mathbb{T}^N$ for later
use.

For $N\geq 1$, we denote by $\{S_n^N\}_{n\geq 0}$ simple random walk
on $\mathbb{T}^N$ such that
\[
P(S_{n+1}^N=y\big|S_n^N=x)=\frac{1}{{\rm deg}(x)}
\]
for each $x\in \mathbb{T}^N$, each neighbor $y$ of $x$ and $n\geq
0$. We assume that $S_0^N=O$. The probability measure and
expectation operator with respect to $\{S_n^N\}_{n\geq 0}$ are
denoted by $\widetilde{P}$ and $\widetilde{E}$.

We define $\Gamma:\mathbb{T}^N\rightarrow Z$ such that $\Gamma(O)=0$
and $\Gamma(y)=\Gamma(x)+1$ when $y$ is a son of $x$. In other
words, for each $x\in \mathbb{T}^N$ there is an unique path $p(O,x)$
from $O$ to $x$ which does not backtrack. $\Gamma(x)$ equals the
length of $p(O,x)$.

\begin{lemma}\label{lemma 4.1}
For any $x\in (0,1]$ and $n\geq 0$,
\[
\widetilde{E}x^{\Gamma(S_n^N)}\leq
\big[\frac{Nx}{N+1}+\frac{1}{(N+1)x}\big]^n.
\]
\end{lemma}

\proof

According to the definition of $S_n^N$,
\begin{align}\label{equ 4.1}
\widetilde{P}\big(\Gamma(S_{n+1}^N)-\Gamma(S_n^N)=1\big|S_n^N=x\big)
=&1-\widetilde{P}\big(\Gamma(S_{n+1}^N)-\Gamma(S_n^N)=-1\big|S_n^N=x\big)\\
=&
\begin{cases}
1 &\text{~if~}x=O,\\
\frac{N}{N+1} &\text{~if~}x\neq O.\\
\end{cases}
\notag
\end{align}
Let $\{Z_n\}_{n\geq 0}$ be a Markov process with state space
$\{\ldots,-2,-1,0,1,2,3,\ldots\}$ and evolve according to
$\{S_n^N\}_{n\geq 1}$. In detail, we assume that $Z_0=0$. For $n\geq
1$, if $S_n^N=O$, then $Z_{n+1}-Z_{n}$ is independent of $S_{n+1}^N$
and satisfies
\[
P(Z_{n+1}-Z_n=1)=1-P(Z_{n+1}-Z_n=-1)=\frac{N}{N+1}.
\]
If $S_n^N\neq O$, then $Z_{n+1}-Z_n=1$ when
$\Gamma(S_{n+1}^N)-\Gamma(S_n^N)=1$ and $Z_{n+1}-Z_n=-1$ when
$\Gamma(S_{n+1}^N)-\Gamma(S_n^N)=-1$.

As a result, for each $n\geq 1$, $Z_{n+1}-Z_n\leq
\Gamma(S_{n+1}^N)-\Gamma(S_n^N)$. Since $Z_0=\Gamma(S_0^N)=0$,
\[
Z_n\leq \Gamma(S_n^N)
\]
for each $n\geq 1$.

Therefore, for $x\in (0,1]$,
\begin{equation}\label{equ 4.2}
\widetilde{E}x^{\Gamma(S_n^N)}\leq Ex^{Z_n}.
\end{equation}
By \eqref{equ 4.1} and the definition of $Z_n$, it is easy to see
that $\{Z_n-Z_{n-1}\}_{n\geq 1}$ are i. i. d random variables such
that
\[
P(Z_n-Z_{n-1}=1)=1-P(Z_n-Z_{n-1}=-1)=\frac{N}{N+1}.
\]
Therefore,
\begin{equation}\label{equ 4.3}
Ex^{Z_n}=\big[Ex^{Z_1-Z_0}\big]^n=\big[\frac{Nx}{N+1}+\frac{1}{(N+1)x}\big]^n.
\end{equation}
Lemma \ref{lemma 4.1} follows from \eqref{equ 4.2} and \eqref{equ
4.3}.

\qed

To control $P(C_t^{O}\neq \emptyset)$ from above, we introduce the
binary contact path process $\{\zeta_t\}_{t\geq 0}$ with random edge
weights on $\mathbb{T}^N$. The classic binary contact path process
is introduced by Griffeath in \cite{Grif1983}, which inspires us a
lot.

The state space of $\{\zeta_t\}_{t\geq 0}$ is
$\{0,1,2,3,\ldots\}^{\mathbb{T}^N}$. At $t=0$, we assume that
$\zeta_0(x)=1$ for each $x\in \mathbb{T}^N$.

When the edge weights $\{\rho(e)\}_{e\in \mathbb{E}^N}$ are given,
$\{\zeta_t\}_{t\geq 0}$ evolves according to Poisson processes
$\{N_x(t):t\geq 0\}_{x\in \mathbb{T}^N}$ and $\{U_{(x,y)}(t):t\geq
0\}_{x\sim y}$. For any $x\in \mathbb{T}^N$, $N_x(\cdot)$ is with
rate $1$. For any $x,y$ such that $x\sim y$, $U_{(x,y)}(\cdot)$ is
with rate $\lambda\rho(x,y)$. Please note that we care the order of
$x$ and $y$, hence $U_{(x,y)}\neq U_{(y,x)}$. We assume that all
these Poisson processes are independent.

For any $t>0$ and $x\in \mathbb{T}^N$, we define
\[
\zeta_{t-}(x)=\lim_{s<t,s\uparrow t}\zeta_s(x)
\]
as the state of $x$ at the moment just before $t$. For $x\in
\mathbb{T}^N$, the state of $x$ may change only at event times of
$N_x(\cdot)$ and $U_{(x,y)}(\cdot)$ for $y\sim x$. At any event time
$s$ of $N_x(\cdot)$, $\zeta_s(x)=0$. At any event time $r$ of
$U_{(x,y)}(\cdot)$, $\zeta_r(x)=\zeta_{r-}(x)+\zeta_{r-}(y)$.

Intuitively, $\{\zeta_t\}_{t\geq 0}$ describes the spread of an
infection disease and the seriousness of the disease for an infected
vertex is considered. An infected vertex $x$ may be further infected
by an infected neighbor $y$. When the infection occurs, the
seriousness of the disease of $y$ will be added to that of $x$.

According to Chapter 9 of \cite{Lig1985}, $\{\zeta_t\}_{t\geq 0}$ is
a linear system with generator $\mathcal{L}$ given by
\begin{align}\label{equ 4.4 generator of binary}
\mathcal{L}f(\zeta)=&\sum_{x\in
\mathbb{T}^N}\big[f(\zeta^{0,x})-f(\zeta)\big]+\sum_{x\in
\mathbb{T}^N}\sum_{y:y\sim x}\lambda
\rho(x,y)\big[f(\zeta^{\zeta(x)+\zeta(y),x})-f(\zeta)\big]
\end{align}
for $f\in C(\{0,1,2,3,\ldots\}^{\mathbb{T}^N})$, where
\[
\zeta^{m,x}(y)=
\begin{cases}
\zeta(y) &\text{~if~} y\neq x,\\
m &\text{~if~} y=x
\end{cases}
\]
for $m\geq 1$ and $x\in \mathbb{T}^N$.

The following lemma is crucial for us to give lower bound of
$\lambda_e$.
\begin{lemma}\label{lemma 4.2}
For any $t\geq 0$,
\[
P_\lambda^N(C_t^{O}\neq \emptyset)\leq E^N_{\lambda}\zeta_t(O).
\]
\end{lemma}

\proof

Let $\{\eta_t\}_{t\geq 0}$ be the contact process defined in
\eqref{equ 1.1 flip rate} with $\eta_0(x)=1$ for any $x\in
\mathbb{T}^N$. Then, according to an approach of graphical
representation introduced in \cite{Har1978}, the contact process
satisfies the dual-relationship that
\begin{equation}\label{equ 4.5}
P_\lambda^{\omega}(C_t^{O}\neq
\emptyset)=P_\lambda^\omega(\eta_t(O)=1)
\end{equation}
for any $\omega\in \Omega$ and therefore
\begin{equation}\label{equ 4.6}
P_\lambda^N(C_t^{O}\neq \emptyset)=P_\lambda^N(\eta_t(O)=1).
\end{equation}
For readers who are not familiar with the self-duality of contact
processes, we give a rigorous proof of \eqref{equ 4.5} in the
appendix.

For any $t\geq 0$ and $x\in \mathbb{T}^N$, we define
\[
\widetilde{\eta}_t(x)=
\begin{cases}
1 &\text{~if~}\zeta_t(x)\geq 1,\\
0&\text{~if~}\zeta_t(x)=0.
\end{cases}
\]
According to the definition of $\{\zeta_t\}_{t\geq 0}$,
$\widetilde{\eta}(x)$ flips from $1$ to $0$ at moment $s$ when and
only when $s$ is an event time of $N_x(\cdot)$ and
$\zeta_{s-}(x)\geq 1$. So $\widetilde{\eta}(x)$ flips from $1$ to
$0$ at rate $1$. $\widetilde{\eta}(x)$ flips from $0$ to $1$ at
moment $r$ when and only when $\zeta_{r-}(x)=0$ and $r$ is an event
time of $U_{(x,y)}(\cdot)$ such that $y\sim x$ and
$\zeta_{r-}(y)\geq 1$. Therefore, $\widetilde{\eta}(x)$ flips from
$0$ to $1$ at rate
\[
\sum_{y:y\sim x}\lambda\rho(x,y)1_{\{\zeta(y)\geq 1\}}=\sum_{y:y\sim
x}\lambda\rho(x,y)\widetilde{\eta}(y).
\]
As a result, $\{\widetilde{\eta}_t\}_{t\geq 0}$ evolves as the same
way as that of $\{\eta_t\}_{t\geq 0}$.

Since $\eta_0(x)=\widetilde{\eta}_0(x)=1$ for each $x\in
\mathbb{T}^N$, $\{\widetilde{\eta}_t\}_{t\geq 0}$ and
$\{\eta_t\}_{t\geq 0}$ have the same probability distribution.

Therefore,
\begin{align}\label{equ 4.7}
P_\lambda^N(\eta_t(O)=1)&=P_\lambda^N(\widetilde{\eta}_t(O)=1)
=P_\lambda^N(\zeta_t(O)\geq 1)\leq E_\lambda^N\zeta_t(O).
\end{align}
Lemma \ref{lemma 4.2} follows from \eqref{equ 4.5} and \eqref{equ
4.7}.

\qed

Finally, we give the proof of $\lambda_e\geq
\big(NE\rho+\frac{M^2}{E\rho}\big)^{-1}$.

\proof[Proof of $\lambda_e\geq
\big(NE\rho+\frac{M^2}{E\rho}\big)^{-1}$]

It is easy to see that we only need to deal with the case where
$M=1$. For general $M>0$, we take $\widetilde{\rho}=\frac{\rho}{M}$
and denote by $\widetilde{\lambda}_e$ the critical value with
respect to $\widetilde{\rho}$. Then,
\[
\lambda_e=\frac{1}{M}\widetilde{\lambda}_e\geq\frac{1}{M}\frac{1}{NE\widetilde{\rho}+\frac{1}{E\widetilde{\rho}}}
=\big(NE\rho+\frac{M^2}{E\rho}\big)^{-1}.
\]
So from now on we assume that $P(\rho\leq 1)=1$.

According to the generator of $\{\zeta_t\}_{t\geq 0}$ given in
\eqref{equ 4.4 generator of binary} and Theorem 9.1.27 of
\cite{Lig1985}, for each $x\in \mathbb{T}^N$ and given edge weights
$\{\rho(e,\omega)\}_{e\in \mathbb{E}^N}$,
\begin{equation}\label{equ 4.8}
\frac{d}{dt}E_\lambda^\omega\zeta_t(x)=-E_\lambda^\omega\zeta_t(x)+\sum_{y:y\sim
x}\lambda\rho(x,y,\omega)E_\lambda^\omega\zeta_t(y).
\end{equation}
For readers who do not want to check the theorem in \cite{Lig1985},
an intuitive explanation of \eqref{equ 4.8} is that $\eqref{equ
4.8}$ is with the form
\[
\frac{d}{dt}Ef(\zeta_t)=E[\mathcal{L}f(\zeta_t)]
\]
with $f(\zeta)=\zeta(x)$ as an `application' of Hille-Yosida
theorem. In fact, Theorem 9.1.27 of \cite{Lig1985} is an extension
of Hille-Yosida theorem to processes of linear systems.

Let $G_\omega$ be $\mathbb{T}^N\times \mathbb{T}^N$ matrix such that
\[
G_\omega(x,y)=
\begin{cases}
\lambda\rho(x,y,\omega) & \text{~if~} x\sim y,\\
0& \text{~otherwise}
\end{cases}
\]
and $I$ be $\mathbb{T}^N\times \mathbb{T}^N$ identity matrix, then
by \eqref{equ 4.8},
\begin{equation}\label{equ 4.9 ODE}
\frac{d}{dt}E_\lambda^\omega\zeta_t=(G_\omega-I)E_\lambda^\omega\zeta_t.
\end{equation}
Since $P(\rho\leq 1)=1$ and there are at most $N+1$ positive
elements in each row of $G_\omega$, it is easy to check that ODE
\eqref{equ 4.9 ODE} satisfies Lipschitz condition under $l_{\infty}$
norm of $R^{\mathbb{T}^N}$ and the series
\[
e^{tG_\omega}=\sum_{n=0}^{+\infty}\frac{t^nG_\omega^n}{n!}
\]
converges. Therefore, according to classic theory of linear ODE, the
unique solution of ODE \eqref{equ 4.9 ODE} is
\begin{equation}\label{equ 4.10}
E_\lambda^\omega\zeta_t=e^{-t}e^{tG_\omega}\zeta_0.
\end{equation}
Since $\zeta_0(x)=1$ for each $x\in \mathbb{T}^N$, by \eqref{equ
4.10},
\begin{equation}\label{equ 4.11}
E_\lambda^\omega\zeta_t(O)=e^{-t}\sum_{n=0}^{+\infty}\sum_{x:x\in
\mathbb{T}^N}\frac{t^nG_\omega^n(O,x)}{n!}.
\end{equation}
For $n\geq 1$, we say that
\[
\overrightarrow{x}=(x_0,x_1,\ldots,x_n)\in
\bigoplus_{j=0}^n\mathbb{T}^N
\]
is a path starting at $O$ with length $n$ when $x_0=O$ and
$x_{j+1}\sim x_j$ for $0\leq j\leq n-1$. Please note that a path may
backtrack.

For $n\geq 1$, we denote by $L_n$ the set of paths starting at $O$
with length $n$.

Then according to the definition of $G_\omega$ and \eqref{equ 4.11},
\begin{equation}\label{equ 4.12}
E_\lambda^\omega\zeta_t(O)=e^{-t}\sum_{n=0}^{+\infty}\frac{t^n\lambda^n}{n!}
\Big(\sum_{\overrightarrow{x}\in
L_n}\prod_{j=0}^{n-1}\rho(x_j,x_{j+1},\omega)\Big),
\end{equation}
where $\overrightarrow{x}=(x_0,x_1,\ldots,x_n)$ and hence
\begin{equation}\label{equ 4.13}
E_\lambda^N\zeta_t(O)=e^{-t}\sum_{n=0}^{+\infty}\frac{t^n\lambda^n}{n!}
\Big(\sum_{\overrightarrow{x}\in
L_n}E\prod_{j=0}^{n-1}\rho(x_j,x_{j+1},\omega)\Big).
\end{equation}

For $\overrightarrow{x}=(x_0,x_1,\ldots,x_n)\in L_n$, there is an
unique path $p(O,x_n)$ from $O$ to $x_n$ with length $\Gamma(x_n)$.
In other words, $p(O,x_n)$ does not backtrack. According to the
structure a tree, the path $\overrightarrow{x}$ contains all the
edges in $p(O,x_n)$. Since $\rho\leq 1$,
\begin{equation}\label{equ 4.14}
E\prod_{j=0}^{n-1}\rho(x_j,x_{j+1},\omega)\leq E\big[\prod_{e\in
p(O,x_n)}\rho(e,\omega)\big]=(E\rho)^{\Gamma(x_n)}.
\end{equation}
Please note that the equation in \eqref{equ 4.14} follows from that
$p(O,x_n)$ is formed with $\Gamma(x_n)$ different edges.

By \eqref{equ 4.13} and \eqref{equ 4.14},
\begin{equation}\label{equ 4.15}
E_\lambda^N\zeta_t(O)\leq
e^{-t}\sum_{n=0}^{+\infty}\frac{t^n\lambda^n}{n!}
\Big[\sum_{\overrightarrow{x}\in L_n}(E\rho)^{\Gamma(x_n)}\Big].
\end{equation}
Since each vertex on $\mathbb{T}^N$ has degree at most $N+1$,
\begin{equation}\label{equ 4.16}
\sum_{\overrightarrow{x}\in L_n}(E\rho)^{\Gamma(x_n)}\leq
(N+1)^n\sum_{\overrightarrow{x}\in
L_n}\prod_{j=0}^{n-1}\frac{1}{{\rm deg}(x_j)}(E\rho)^{\Gamma(x_n)}.
\end{equation}
By the definition of $\{S_n^N\}_{n\geq 1}$, for
$\overrightarrow{x}=(x_0,x_1,\ldots,x_n)\in L_n$,
\begin{equation*}
\widetilde{P}(S_j^N=x_j,0\leq j\leq
n)=\prod_{j=0}^{n-1}\frac{1}{{\rm deg}(x_j)}
\end{equation*}
and hence
\begin{equation}\label{equ 4.17}
\sum_{\overrightarrow{x}\in L_n}\prod_{j=0}^{n-1}\frac{1}{{\rm
deg}(x_j)}(E\rho)^{\Gamma(x_n)}=\widetilde{E}\big[(E\rho)^{\Gamma(S_n)}\big].
\end{equation}

By \eqref{equ 4.16} and \eqref{equ 4.17},
\begin{equation}\label{equ 4.18}
\sum_{\overrightarrow{x}\in L_n}(E\rho)^{\Gamma(x_n)}\leq
(N+1)^n\widetilde{E}\big[(E\rho)^{\Gamma(S_n)}\big].
\end{equation}
By \eqref{equ 4.15} and \eqref{equ 4.18},
\begin{equation}\label{equ 4.19}
E_\lambda^N\zeta_t(O)\leq
e^{-t}\sum_{n=0}^{+\infty}\frac{t^n\lambda^n(N+1)^n}{n!}\widetilde{E}\big[(E\rho)^{\Gamma(S_n)}\big].
\end{equation}
By \eqref{equ 4.19} and Lemma \ref{lemma 4.1},
\begin{align}\label{equ 4.20}
E_\lambda^N\zeta_t(O)&\leq e^{-t}\sum_{n=0}^{+\infty}\frac{t^n\lambda^n(N+1)^n}{n!}\big[\frac{NE\rho}{N+1}+\frac{1}{(N+1)E\rho}\big]^n\notag\\
&=\exp\Big\{t\big[\lambda(NE\rho+\frac{1}{E\rho})-1\big]\Big\}.
\end{align}
By Lemma \ref{lemma 4.2} and \eqref{equ 4.20},
\[
P_\lambda^N(C_t^{O}\neq \emptyset)\leq
\exp\Big\{t\big[\lambda(NE\rho+\frac{1}{E\rho})-1\big]\Big\}.
\]
Therefore,
\[
\limsup_{t\rightarrow+\infty}\frac{1}{t}\log P_\lambda^N(C_t^{O}\neq
\emptyset)\leq \lambda(NE\rho+\frac{1}{E\rho})-1<0
\]
when
\[
\lambda<(NE\rho+\frac{1}{E\rho})^{-1}.
\]
As a result,
\[
\lambda_e\geq (NE\rho+\frac{1}{E\rho})^{-1}.
\]

\qed

Now we can complete the proof Theorem \ref{theorem 2.0} and Theorem
\ref{theorem 2.1 limit theorem for two critical}.

\proof[Proof of Theorem \ref{theorem 2.0}]

According to Corollary \ref{corollary 3.2}, we only need to show
that $\lambda_c>0$ in any case. Since
\[
\lambda_c\geq \lambda_e\geq (NE\rho+\frac{M^2}{E\rho})^{-1}>0,
\]
the proof is complete.

\qed

\proof[Proof of Theorem \ref{theorem 2.1 limit theorem for two
critical}]

Since $\lambda_e\geq (NE\rho+\frac{M^2}{E\rho})^{-1}$,
\[
\liminf_{N\rightarrow+\infty}N\lambda_e(N)\geq \frac{1}{E\rho}.
\]
Since $\lambda_e\leq \lambda_c$ and we have shown that
\[
\limsup_{N\rightarrow+\infty}N\lambda_c(N)\leq \frac{1}{E\rho}
\]
in Section \ref{section three},
\[
\lim_{N\rightarrow+\infty}N\lambda_e(N)=\lim_{N\rightarrow+\infty}N\lambda_c(N)=\frac{1}{E\rho}
\]
and the proof is complete.

\qed

\section{Critical value under quenched measure}\label{section five}
In this section we discuss the critical value under quenched
measure. For later use, we identify $\mathbb{T}^N$ with the set
\[
\{O\}\bigcup\bigcup_{m=1}^{+\infty}\{1,2,3,\ldots,N\}^m.
\]
In detail, $O$ is the root of $\mathbb{T}^N$. For $1\leq j\leq N$,
$j$ represents the $j$th son of $O$. For $m\geq 1$, $1\leq j\leq N$
and
\[
(k_1,k_2,\ldots,k_m)\in\{1,2,\ldots,N\}^m,
\]
$(k_1,k_2,\ldots,k_m,j)$ represents the $j$th son of
$(k_1,k_2,\ldots,k_m)$. The following picture describes the first
three generations of $\mathbb{T}^2$.

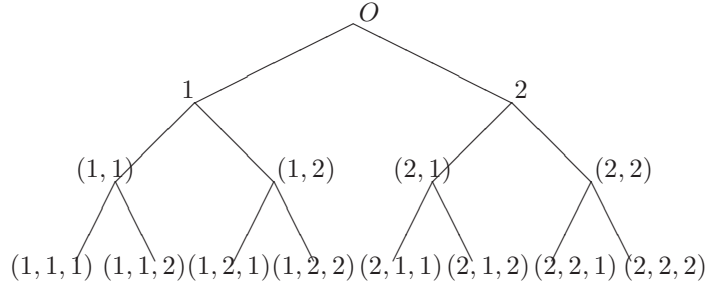
\begin{center}
\begin{picture}(250,110)
\put(112,96){\(O\)} \put(110,95){\line(2,-1){60}}
\put(110,95){\line(-2,-1){60}} \put(170,65){\line(1,-1){30}}
\put(170,65){\line(-1,-1){30}} \put(50,65){\line(1,-1){30}}
\put(50,65){\line(-1,-1){30}} \put(20,35){\line(1,-2){15}}
\put(20,35){\line(-1,-2){15}} \put(80,35){\line(1,-2){15}}
\put(80,35){\line(-1,-2){15}} \put(140,35){\line(1,-2){15}}
\put(140,35){\line(-1,-2){15}} \put(200,35){\line(1,-2){15}}
\put(200,35){\line(-1,-2){15}} \put(45,67){$1$} \put(171,67){$2$}
\put(5,37){$(1,1)$} \put(81,37){$(1,2)$} \put(125,37){$(2,1)$}
\put(201,37){$(2,2)$} \put(-20,0){$(1,1,1)$} \put(15,0){$(1,1,2)$}
\put(47,0){$(1,2,1)$} \put(79,0){$(1,2,2)$} \put(112,0){$(2,1,1)$}
\put(145,0){$(2,1,2)$} \put(178,0){$(2,2,1)$} \put(212,0){$(2,2,2)$}
\end{picture}
\end{center}

For each $1\leq j\leq N$, we define injection
$\varphi_j:\mathbb{T}^N\rightarrow \mathbb{T}^N$ such that
\[
\varphi_j(O)=j
\]
and
\[
\varphi_j(k_1,k_2,\ldots,k_m)=(j,k_1,k_2,\ldots,k_m)
\]
for each $m\geq 1$ and any $(k_1,k_2,\ldots,k_m)\in
\{1,2,\ldots,N\}^m$.

For $e\in \mathbb{E}^N$ with endpoints $x,y\in \mathbb{T}^N$, we
denote by $e_j$ the edge with endpoints $\varphi_j(x)$ and
$\varphi_j(y)$. For $\omega\in \Omega_N$ and $j\geq 1$, we denote by
$\omega_j$ the sample point such that
\[
\rho(e,\omega_j)=\rho(e_j,\omega)
\]
for each $e\in \mathbb{E}^N$. That is to say, if $\mathbb{T}^N$ is
with edge weights $\{\rho(e,\omega)\}_{e\in \mathbb{E}^N}$, then $j$
and its descendants form a regular tree which is rooted at $j$ and
with edge weights $\{\rho(e,\omega_j)\}_{e\in \mathbb{E}^N}$.

For any $\lambda>0$, $N\geq 1$ and $1\leq j\leq N$, we define
\[
H(\lambda,N)=\{\omega\in \Omega_N:P_\lambda^\omega(\forall~t\geq
0,C_t^{O}\neq\emptyset)=0\}
\]
and
\[
H(\lambda,N,j)=\{\omega\in
\Omega_N:P_\lambda^{\omega_j}(\forall~t\geq
0,C_t^{O}\neq\emptyset)=0\}.
\]
The following lemma shows that $H(\lambda,N)$ satisfies a zero-one
law, which is crucial for us to prove Theorem \ref{theorem 2.2 two
critical values are equal}. Please note that $A_N$ in the lemma is
the same as that defined in \eqref{equ 2.8}.

\begin{lemma}\label{lemma 5.1}
If $P(\rho>0)<1$ and $N>1/P(\rho>0)$, then $0<\mu_N(A_N)<1$ and
\[
\mu_N\big(H(\lambda,N)\big)\in \{1-\mu_N(A_N),1\}
\]
for any $\lambda>0$.

If $P(\rho>0)=1$ and $N\geq 2$, then
\[
\mu_N\big(H(\lambda,N)\big)\in \{0,1\}
\]
for any $\lambda>0$.
\end{lemma}

\proof

For any $\omega\in \Omega$, we define
\[
B(\omega)=\{1\leq j\leq N:\rho(O,j,\omega)>0\}
\]
as the set of sons which $O$ can infect.

According to the strong Markov property, for $1\leq j\leq N$,
\begin{equation}\label{equ 5.1}
P_\lambda^\omega(\forall~t\geq 0,C_t^{O}\neq\emptyset)\geq
P_\lambda^\omega(\exists~t>0, j\in
C_t^{O})P_\lambda^{\omega_j}(\forall~t\geq 0,C_t^{O}\neq\emptyset).
\end{equation}
If $j\in B(\omega)$, then $P_\lambda^{\omega}(\exists~t>0, j\in
C_t^{O})>0$. Therefore, by \eqref{equ 5.1},
$P_\lambda^\omega(\forall~t\geq 0,C_t^{O}\neq\emptyset)=0$ and $j\in
B(\omega)$ implies that $P_\lambda^{\omega_j}(\forall~t\geq
0,C_t^{O}\neq\emptyset)=0$. As a result,
\begin{equation}\label{equ 5.2}
H(\lambda,N)\subseteq \{\omega:\omega\in \bigcap_{j\in
B(\omega)}H(\lambda,N,j)\}.
\end{equation}

Since $\{\rho(e)\}_{e\in \mathbb{E}^N}$ are i.i.d,
$H(\lambda,N,1),H(\lambda,N,2),\ldots,H(\lambda,N,N)$ are
independent of $B(\omega)$ and are i.i.d events which have the same
probability distribution as that of $H(\lambda,N)$ under $\mu_N$.
Therefore, by \eqref{equ 5.2},
\begin{equation}\label{equ 5.3}
\mu_N\big(H(\lambda,N)\big)\leq
\sum_{k=0}^{N}p_k\Big[\mu_N\big(H(\lambda,N)\big)\Big]^k,
\end{equation}
where
\begin{align*}
p_k&=\mu_N(\omega:|B(\omega)|=k)\\
&={N\choose k}P(\rho>0)^k\big(1-P(\rho>0)\big)^{N-k}.
\end{align*}
For $x\in [0,1]$, we define
\[
f(x)=\sum_{k=0}^{N}p_kx^k.
\]
As we have shown in Section \ref{section two}, $D(\omega)$ defined
in \eqref{equ 2.7 two} is a Galton-Watson tree with binomial
offspring distribution $B(N,P(\rho>0))$ and $1-\mu_N(A_N)$ is the
extinction probability of $D(\omega)$.

When $P(\rho>0)<1$ and $N>1/P(\rho>0)$, the mean of $B(N,P(\rho>0))$
is larger than one. Then according to the extinction criterion of
Galton-Watson trees, $1-\mu_N(A_N)$ is the unique solution in
$(0,1)$ to the equation $x=f(x)$ and $f(y)<y$ for $y\in
\big(1-\mu_N(A_N),1\big)$. By \eqref{equ 5.3},
$\mu_N\big(H(\lambda,N)\big)\leq
f\Big(\mu_N\big(H(\lambda,N)\big)\Big)$, hence
\begin{equation}\label{equ 5.4}
\mu_N\big(H(\lambda,N)\big)\in [0,1-\mu_N(A_N)]\cup\{1\}.
\end{equation}
For any $\omega\in \Omega_N\setminus A_N$, $|D(\omega)|$ is finite
and hence the Markov process $\{C_t^{O}\}_{t\geq 0}$ under the
measure $P_\lambda^\omega$ is with finite state space
$\{A:A\subseteq D(\omega)\}$ and unique absorption state
$\emptyset$, which makes $\{C_t^{O}\}_{t\geq 0}$ frozen in
$\emptyset$ eventually. As a result,
\[
P_\lambda^\omega(\forall~t\geq 0,C_t^{O}\neq\emptyset)=0
\]
for any $\omega\in \Omega_N\setminus A_N$ and hence
\begin{equation}\label{equ 5.5}
\mu_N\big(H(\lambda,N)\big)\geq \mu_N(\Omega_N\setminus
A_N)=1-\mu_N(A_N).
\end{equation}
By \eqref{equ 5.4} and \eqref{equ 5.5},
\[
\mu_N\big(H(\lambda,N)\big)\in \{1-\mu_N(A_N),1\}.
\]

When $P(\rho>0)=1$ and $N\geq 2$, \eqref{equ 5.3} turns into
\begin{equation}\label{equ 5.6}
\mu_N\big(H(\lambda,N)\big)\leq
\Big[\mu_N\big(H(\lambda,N)\big)\Big]^N.
\end{equation}
If $0<\mu_N\big(H(\lambda,N)\big)<1$, then
\[
\Big[\mu_N\big(H(\lambda,N)\big)\Big]^N<\mu_N\big(H(\lambda,N)\big)
\]
since $N\geq 2$, which is contradictory to \eqref{equ 5.6}.
Therefore,
\[
\mu_N\big(H(\lambda,N)\big)\in \{0,1\}.
\]

\qed

In the case where $P(\rho>0)=1$ and $N=1$, \eqref{equ 5.3} turns
into $\mu\big(H(\lambda,N)\big)\leq \mu\big(H(\lambda,N)\big)$,
which gives no information. This is why this case should be
discussed specially. We propose an open question about the critical
value in this case in section \ref{section six}.

Finally we give the proof of Theorem \ref{theorem 2.2 two critical
values are equal}.

\proof[Proof of Theorem \ref{theorem 2.2 two critical values are
equal}]

We first consider the case where $P(\rho>0)=1$ and $N\geq 2$. In
this case, we have shown in the proof of Theorem \ref{theorem 2.0}
that
\[
\lambda_c(N)\in (0,+\infty).
\]
So we only need to show that
$\widehat{\lambda}_c(\omega,N)=\lambda_c(N)$ with probability one.
For $m>1/\lambda_c(N)$, let $\lambda_m=\lambda_c(N)-\frac{1}{m}$ and
$\beta_m=\lambda_c(N)+\frac{1}{m}$, then according to the definition
of $\lambda_c(N)$,
\[
P_{\lambda_m}^N(\forall~t\geq0,C_t^{O}\neq
\emptyset)=E_{\lambda_m}^N\big[P_{\lambda_m}^{\omega}(\forall~t\geq0,C_t^{O}\neq\emptyset)\big]=0
\]
and
\[
P_{\beta_m}^N(\forall~t\geq0,C_t^{O}\neq
\emptyset)=E_{\beta_m}^N\big[P_{\beta_m}^{\omega}(\forall~t\geq0,C_t^{O}\neq\emptyset)\big]>0.
\]
Therefore, according to lemma \ref{lemma 5.1},
\begin{equation}\label{equ 5.7}
\mu_N\big(H(\lambda_m,N)\big)=1
\end{equation}
and
\begin{equation}\label{equ 5.8}
\mu_N\big(H(\beta_m,N)\big)=0.
\end{equation}
Let
\[
K_N=\bigcap_{m}H(\lambda_m,N)\bigcap\bigcap_{m}\big(\Omega_N\setminus
H(\beta_m,N)\big),
\]
then
\[
\mu_N(K_N)=1
\]
according to \eqref{equ 5.7} and \eqref{equ 5.8}. For $\omega\in
K_N$,
\[
P_{\lambda_m}^{\omega}(\forall~t\geq0,C_t^{O}\neq\emptyset)=0,
P_{\beta_m}^{\omega}(\forall~t\geq0,C_t^{O}\neq\emptyset)>0
\]
and hence
\[
\lambda_m\leq\widehat{\lambda}_c(\omega,N)\leq \beta_m.
\]
Let $m\rightarrow+\infty$, then we have that
\[
\widehat{\lambda}_c(\omega,N)=\lambda_c(N)
\]
for $\omega\in K_N$.

Now we deal with the case where $P(\rho>0)<1$ and $N>1/P(\rho>0)$.
As we have shown in the proof of Theorem \ref{theorem 2.0},
\[
\lambda_c(N)\in (0,+\infty)
\]
in this case. We also use $\lambda_m$ and $\beta_m$ to denote
$\lambda_c(N)-\frac{1}{m}$ and $\lambda_c(N)+\frac{1}{m}$
respectively. According to a similar analysis with that of the first
case and Lemma \ref{lemma 5.1},
\begin{equation}\label{equ 5.9}
\mu_N\big(H(\lambda_m,N)\big)=1
\end{equation}
and
\begin{equation}\label{equ 5.10}
\mu_N\big(H(\beta_m,N)\big)=1-\mu_N(A_N).
\end{equation}
We have shown in the proof of Lemma \ref{lemma 5.1} that
\[
H(\lambda,N)\supseteq \Omega_N\setminus A_N
\]
for any $\lambda>0$, hence by \eqref{equ 5.10},
\begin{equation}\label{equ 5.11}
\mu_N\big(H(\beta_m,N)\cap A_N\big)=0.
\end{equation}
Let
\[
K_N=\Big(A_N\setminus\bigcup_mH(\beta_m,N)\Big)\bigcap\bigcap_mH(\lambda,N),
\]
then $K_N\subseteq A_N$ and $\mu_N(A_N\setminus K_N)=0$ according to
\eqref{equ 5.9} and \eqref{equ 5.11}.

According to a similar analysis with that of the first case, it is
easy to see that
\[
\widehat{\lambda}_c(\omega,N)=\lambda_c(N)
\]
for any $\omega\in K_N$.

For $\omega\in \Omega_N\setminus A_N$, $|D(\omega)|<+\infty$ and
hence
\[
P_{\lambda}^{\omega}(\forall~t\geq0,C_t^{O}\neq\emptyset)=0
\]
for any $\lambda>0$ as we have shown in the proof of Lemma
\ref{lemma 5.1}. As a result.
\[
\widehat{\lambda}_c(\omega,N)=+\infty
\]
for any $\omega\in \Omega_N\setminus A_N$.

For the last case where $P(\rho>0)<1$ and $N\leq 1/P(\rho>0)$,
\[
A_N=\emptyset
\]
according to the extinction criterion of Galton-Watson trees, and
hence
\[
\widehat{\lambda}_c(\omega,N)=+\infty
\]
for any $\omega\in \Omega_N$.

\qed

\section{An Open question for $N=1$}\label{section six}
When $N=1$, our model turns into the contact process with random
edge weights on $\mathbb{Z}$. We do not manage to give an criterion
to judge whether $\lambda_c<+\infty$ in this case.

There are two trivial cases for this problem. If $P(\rho>0)<1$, then
$|D(\omega)|$ is finite with probability one and hence
$\lambda_c=+\infty$. If $P(\rho>\epsilon)=1$ for some $\epsilon>0$,
then $\lambda_c<+\infty$ since the classic contact process on
$\mathbb{Z}$ has finite critical value (see \cite{Lig1995}). So we
only need to deal with the case where $P(\rho>0)=1$ but
$P(\rho<x)>0$ for any $x\in (0,1)$.

We do not think that $P(\rho>0)=1$ is sufficient or
$P(\rho>\epsilon)=1$ is necessary for finite critical value. We
guess that the probability of $\{\rho<x\}$ for small $x$ is crucial.

To make our question concrete, for $\alpha>0$, we assume that
\begin{equation}\label{equ 6.1}
P(\rho<x)=x^\alpha
\end{equation}
for $x\in (0,1)$. We denote by $\lambda_c(\alpha)$ the critical
value with respect to $\rho$. Then it is obviously that
\[
\lambda_c(\alpha_1)\leq \lambda_c(\alpha_2)
\]
for $\alpha_1>\alpha_2$.

Then it is reasonable to ask the following question.

\begin{question}\label{question 6.1}
We assume that $N=1$ and $\rho$ has the distribution as that in
\eqref{equ 6.1}. Then is there a critical value $0<\alpha_c<+\infty$
such that
\[
\lambda_c(\alpha)<+\infty
\]
for $\alpha>\alpha_c$ and
\[
\lambda_c(\alpha)=+\infty
\]
for $\alpha<\alpha_c$?
\end{question}

If the answer to Question \ref{question 6.1} is positive, a further
problem is how to estimate $\alpha_c$, which will bring more
interesting work to do. We will work on Question \ref{question 6.1}
as a further study and hope to discuss with readers who are
interested in this question.

\appendix{}
\section{Appendix}
\proof [Proof of \eqref{equ 4.5}] According to the flip rate
functions given by \eqref{equ 1.1 flip rate}, the Markov process
$\{C_t\}_{t\geq 0}$ is with state space
$2^{\mathbb{T}^N}:=\{A:A\subseteq \mathbb{T}^N\}$ and has generator
given by
\begin{equation}\label{equ A.1}
\mathcal{L}f(A)=\sum_{x:x\in A}\big[f(A\setminus x)-f(A)\big]
+\sum_{x:x\in A}\sum_{y:y\sim
x}\lambda\rho(x,y)\big[f(A\cup\{y\})-f(A)\big]
\end{equation}
for $f\in C(2^{\mathbb{T}^N})$ and $A\subseteq \mathbb{T}^N$.

We define $H:2^{\mathbb{T}^N}\times
2^{\mathbb{T}^N}\rightarrow\{0,1\}$ that
\begin{equation}\label{equ A.2}
H(A,B)=
\begin{cases}
1 &\text{~if~} A\cap B=\emptyset,\\
0 &\text{~if~} A\cap B\neq \emptyset
\end{cases}
\end{equation}
for $A,B\subseteq \mathbb{T}^N$.

By \eqref{equ A.2},
\begin{equation}\label{equ A.3}
H(A,B\cup C)=H(A,B)H(A,C)
\end{equation}
for $A,B,C\subseteq T^N$.

By \eqref{equ A.1}, \eqref{equ A.3} and direct calculation,
\begin{align}\label{equ A.4}
\mathcal{L}H(A,\cdot)(B)=&\sum_{x:x\in B}\big[H(A,B\setminus
x)-H(A,B)\big] \notag\\
&+\sum_{x:x\in B}\sum_{y:y\sim
x}\lambda\rho(x,y)\big[H(A,B\cup\{y\})-H(A,B)\big]\notag\\
=&\sum_{x:x\in B}H(A,B\setminus x)\big[1-H(A,\{x\})\big]\notag\\
&+\sum_{x:x\in B}\sum_{y:y\sim
x}\lambda\rho(x,y)H(A,B)\big[H(A,\{y\})-1\big]\notag\\
=&\sum_{x:x\in A\cap B}H(A,B\setminus x)-\sum_{x:x\in
B}\sum_{y:y\sim x,\atop y\in A}\lambda\rho(x,y)H(A,B)
\end{align}
for $A,B\subseteq \mathbb{T}^N$. According to a similar calculation,
\begin{align}\label{equ A.5}
\mathcal{L}H(\cdot,B)(A)&=\sum_{x:x\in A\cap B}H(A\setminus
x,B)-\sum_{y:y\in A}\sum_{x:x\sim y,\atop x\in
B}\lambda\rho(x,y)H(A,B)\notag\\
&=\sum_{x:x\in A\cap B}H(A\setminus x,B)-\sum_{x:x\in
B}\sum_{y:y\sim x,\atop y\in A}\lambda\rho(x,y)H(A,B).
\end{align}

It is easy to see that
\[
H(A,B\setminus x)=H(A\setminus x,B)
\]
for $x\in A\cap B$. Therefore, by \eqref{equ A.4} and \eqref{equ
A.5},
\begin{equation}\label{equ A.6}
\mathcal{L}H(\cdot,B)(A)=\mathcal{L}H(A,\cdot)(B)
\end{equation}
for $A,B\subseteq \mathbb{T}^N$.

We write $C_t$ as $C_t^A$ when $C_0=A$. Then, according to
\eqref{equ A.6} and Theorem 3.39 of \cite{Lig2010},
\begin{equation}\label{equ A.7}
E_\lambda^\omega H(A,C_t^{B})=E_\lambda^\omega H(C_t^{A},B)
\end{equation}
for $A,B\subseteq \mathbb{T}^N$ and $t\geq 0$. Let $A=\{O\}$ and
$B=\mathbb{T}^N$, then \eqref{equ 4.5} follows from \eqref{equ A.7}.

\qed

\textbf{Acknowledgments.} The author is grateful to the financial
support from the National Natural Science Foundation of China with
grant number 11171342 and China Postdoctoral Science Foundation (No.
2015M571095).

{}
\end{document}